\documentclass[11pt]{article}
\usepackage{amsmath,amssymb}
\usepackage{amsthm}
\usepackage{epsf}
\usepackage{graphicx}

\newtheorem{lem}{Lemma}[section]
\newtheorem{thm}{Theorem}[section]

\theoremstyle{remark}

\theoremstyle{definition}

\numberwithin{equation}{section}

\begin{document}

\title{Exponential and Algebraical Stability of Traveling Wavefronts in Periodic
Spatial-Temporal Environments}
\author{Ming Mei\thanks{%
Department of Mathematics, Champlain College Saint-Lambert, Saint-Lambert,
Quebec J4P 3P2, Canada, and Department of Mathematics and Statistics, McGill
University, Montreal, Quebec H3A 2K6, Canada, \texttt{mei@math.mcgill.ca}.
This work is supported in part by the NSERC of Canada. },\,\, Chunhua Ou%
\thanks{%
Department of Mathematics and Statistics, Memorial University of
Newfoundland, St. John's, NL, Canada A1C 5S7, \texttt{ou@mun.ca}. This work
is supported in part by the NSERC grant of Canada. The corresponding author.}%
\, \, and \, Xiao-Qiang Zhao\thanks{
Department of Mathematics and Statistics, Memorial University of
Newfoundland, St. John's, NL, Canada A1C 5S7, \texttt{zhao@mun.ca}. This
work is supported in part by the NSERC grant of Canada and the MITACS of
Canada.}}
\date{}
\maketitle

\begin{abstract}
Global stability of traveling wavefronts in a periodic spatial-temporal
environment in $n$-dimension ($n\ge 1$) is studied.   The wavefront is proved to be exponentially stable in the form of $%
O(e^{-\mu t})$ for some $\mu>0$, when the wave speed is greater than
the critical one, and algebraically stable in the form of
$O(t^{-n/2})$ in the critical case. A new and easy to follow method
is developed. These results are then extended to the case of
time-periodic media. Finally, we illustrate how the stability result
can be directly used to obtain the uniqueness of the wavefront with
a given speed.
\end{abstract}



\smallskip


\smallskip

\textbf{AMS 2000 MSC}: 35K57, 92D25 \baselineskip=18pt

\section{Introduction}

Since the pioneering works of Fisher \cite{F} and KPP \cite{KPP},
wave propagations for reaction-diffusion models in fluid mechanics,
physical, chemical and biological systems have been investigated
intensively. For example, such wave behaviors with front propagation
into unstable states can be frequently seen in fluid dynamical
experiments on Taylor-Couette\cite {AC,LMW} and Rayleigh-Benard
system \cite{FS}, and in chemical wave experiment \cite{SS,VRCY} as
well as in population dynamics, combustion and biological invasions
in a \textit{homogeneous} environment \cite{AW1,AW2}. The theory of
spreading speeds and monostable traveling waves has been developed
to monotone semiflows in such a way that it can be applied to
various evolution systems admitting the comparison principle (see
\cite {Wein1982,Lui,LiangZhao}). For the stability of monostable
waves including the minimal wave fronts in some regular and
time-delayed reaction-diffusion equations in spatially homogeneous
habitat, we refer to \cite{VVV,gourley,M1,M2,M3,MOZ} and the
references therein. The stability for bistable waves to the regular
or time-delayed reaction-diffusion equations has also been
intensively studied by many people, for example, see
\cite{Mac-F,Matano,Kapitula,Smith-Zhao,X,VVV} and the references
therein.

Recently, the study of traveling waves and spreading speeds in \textit{%
heterogeneous} media or periodically fragmented environments has been
drawing more attention, see, e.g., the survey paper \cite{X}. In a series of
works, Berestyski et al. \cite{B2005,B2008,B2010} and Hamel \cite{H2008}
established the existence of traveling fronts of the following equation
\begin{equation}
u_t-\nabla \cdot (A(x)\nabla u)=f(x,u),\;\;\;x\in \mathbb{R}^n,
\label{1}
\end{equation}
where the diffusion matrix $A(x)$ and the reaction term $f$ both are
periodic in $x=(x_1,x_2,...,x_n)$ in the sense that there exist constants $%
L_i>0,\;1\leq i\leq n$,  such that
\[
A(x+L_i{\bf e}_i)=A(x), \quad f(x+L_i{\bf e}_i,u)=f(x,u), \quad
\forall 1\leq i\leq n,
\]
with ${\bf e}_i$ being the $i$-th vector in the standard base of
$\mathbb{R}^n$ (i.e., its $i$-th component is $1$ and all other
components are $0$). The matrix $A(x)$ is also assumed to be
symmetric
and lies in the space $C^{1,\alpha }(\mathbb{R}^n).$ The reaction term $%
f(x,u)$ is a continuous function of class $C^{0,\alpha }$ with respect to $x$
locally uniformly in $u\in \mathbb{R}$. The derivative $\frac{\partial f}{%
\partial u}$ is assumed to exist and be continuous in $u \in \mathbb{R}$ uniformly for $x$ in $\mathbb{R}%
^n.$ For instance, $f$ can be taken as $f(x,u)=u(\mu (x)-u), $ with
the positive periodic function $\mu (x)$ usually being the
carrying capacity of the periodic environment. Furthermore, for the matrix $%
A,$ it is usually assumed to be bounded, and uniformly elliptic in the sense
that
\begin{equation}
\sum_{1\leq i,j\leq n}a_{i,j}(x)\xi _i\xi _j\geq \alpha _0|\xi
|^2,\;\;\alpha _0>0  \label{1a}
\end{equation}
for any vector $\xi \in \mathbb{R}^n.$

In the special case where the period $L_i=0,$ namely, the matrix $A$
and the reaction term $f$ are independent of $x,$ model (\ref{1}) is
a generalization of the simple Fisher-KPP model
\begin{equation}
\left\{
\begin{array}{l}
u_t-\Delta u=f(u)\vspace{3.01pt} \\
f(0)=f(1)=0,\;\;f(u)>0,\,\forall u\in (0,1).
\end{array}
\right.  \label{2}
\end{equation}
It is well known that for any $c\geq 2$, system (\ref{2}) with
$f(u)=u(1-u)$ has a planar traveling wave solution $u(t,x)=U(x\cdot
\mathbf{e}+ct)$ with $|\mathbf{e}|=1$, where $-\mathbf{e}$ is the
direction of propagation, while it has no such traveling wave for
any $c<2$.

The existence of traveling fronts to (\ref{1})  has been studied
recently in \cite{B2005,H2008}. For spreading speeds and traveling waves of
other types of biological evolution systems in a periodic habitat, we
further refer to \cite{Wein2002,JinZhao, WengZhao} and the references
therein. However, as mentioned in the survey paper \cite[page 185]{X}, the stability analysis
of traveling wavefronts in periodic media is a very challenging research
topic. In particular, the author of \cite{X} emphasized that the dynamics of
the slow waves moving at the minimal speeds requires ``a more delicate
argument and awaits further investigation''.

In this paper, we will concentrate on the stability and uniqueness
of the wavefront of (\ref{1}), an extension of the classical
KPP-Fisher model in heterogeneous media.  Motivated by the work of
\cite{Moet} on the classical Fisher-KPP equation and our recent
study \cite{MOZ} on nonlocal time-delayed reaction-diffusion
equations, we further investigate the stability of wavefronts for
the periodic spatial-temporal reaction-diffusion equations by the
Green function method. A \emph{new and easy to follow} approach in
terms of a wave transformation is developed. The wavefront is proved
to be exponentially stable in the case where the wave-speed is
greater than the minimal one (also called the critical speed), and
algebraically stable in the case of the minimal speed. Furthermore,
we extend our results to the case of  time-periodic media. In the
final concluding remark, we illustrate how the stability result can
be directly used to obtain the uniqueness of the wavefront with a
given speed.

After we finished our first draft of this paper, we happened to
obtain an important  preprint \cite{H2} from Drs. Roques and Hamel's
personal websites, where
they investigated the stability and uniqueness of pulsating wavefronts for (%
\ref{1}) as well. It is seen that both papers are focusing on the
open problem of the stability of wavefronts.  However, there exist
some major differences.
Firstly, our developed method is quite different from
that in \cite{H2}. Secondly, we prove that all non-critical
pulsating wavefronts are exponentially stable in the form of
$O(e^{-\mu t})$ for some $\mu
>0$, and the critical pulsating traveling wave is algebraically
stable in the form of $O(t^{-n/2})$, while there are no convergence
rates in the global stability results of \cite{H2}. Thirdly, we
assume that the product of the initial perturbation and the weight
function is in $L^1(R^n)$, while it is assumed in \cite{H2} that
initial functions behave like a front.

\section{Exponential and algebraical stability}

In this section, we are going to investigate the global stability of planar
traveling waves to the spatial periodic equation \eqref{1}. Let $%
w(x,x\cdot \mathbf{e}+ct)$ be a given planar traveling wavefront
with wave speed $c$ for some direction $\mathbf{e}$
($|\mathbf{e}|=1,$ $-\mathbf{e}$ is the direction of propagation)
such that $w(x+L_i{\bf e}_i,\xi)=w(x,\xi), \, \, \forall x\in
\mathbb{R}^n,\,  \xi=x\cdot \mathbf{e}+ct \in \mathbb{R},\, 1\leq
i\leq n$. Without loss of generality, we can always assume that
$\mathbf{e}=\mathbf{e}_1=(1,0,...,0)$ after a rotation of
coordinates. Note that the existence and the asymptotic behavior of
$w(x,\xi ),$ $\xi =x_1+ct,$ can be obtained from \cite
{B2005,H2008}.

For equation (\ref{1}) with $f(x,0)=0$ to have wavefronts connecting
two steady states, the following sub-linearity of the function $f$
is also assumed in \cite{B2005,H2008}:

\begin{equation}
s\rightarrow \frac{f(x,s)}s\text{ is decreasing in }s>0\text{ for
any }x\in \mathbb{R}^n\mathbf{,}  \label{2.1}
\end{equation}
and
\begin{equation}
f(x,s)\leq 0,\;s\geq M\text{ for some positive constant }M.  \label{2.2}
\end{equation}

We have the following simple observation.

\begin{lem}
If condition (\ref{2.1}) holds, then we have
\begin{equation}\label{2.1a}
f_u(x,u)\leq f_u(x,0),\;\;\;\forall x\in \mathbb{R}^n,\, u\geq 0,
\end{equation}
where $f_u(x,u)$ is the partial derivative of $f$ with respect to
the variable $u$.
\end{lem}

\begin{proof}
Define
\[
\quad g(x,s):=\frac{f(x,s)}s,\, \, \forall x\in \mathbb{R}^n, \,
s>0.
\]
Then  \eqref{2.1} implies that the function $g$ is decreasing in
$s>0$, and hence, $\frac{\partial g}{\partial s}\leq 0$. Since
$f(x,s)=sg(x,s)$, we have
\[
\frac{\partial f}{\partial s}=g(x,s)+s\frac{\partial g}{\partial
s}\leq g(x,s)\leq g(x,0):=\frac{\partial f(x,0)}{\partial s}.
\]
This completes the proof.
\end{proof}

The non-negative steady-states of (\ref{1}), involving in the
boundary condition of wavefront at $\pm \infty ,$ are the
$C^{2,\alpha }(\mathbb{R}^n)$-solutions $p(x)$ defined by
\[
\left\{
\begin{array}{l}
\ \nabla \cdot (A(x)\nabla p(x))+f(x,p(x))=0,\vspace{3pt} \\
p(x)\geq 0,\;\;p(x+L_i{\bf e}_i)=p(x), \, \forall 1\leq i\leq n.
\end{array}
\right.
\]
Obviously, the constant $0$ is a trivial steady-state and its local
stability is associated with an eigenvalue problem. Let $\mu _1$ be the
principal eigenvalue of the operator $L_0$ defined by
\[
L_0\phi =\ \nabla \cdot (A(x)\nabla \phi )+f_u(x,0)\phi ,
\]
that is,
\begin{equation}
\left\{
\begin{array}{l}
\ \nabla \cdot (A(x)\nabla \phi )+f_u(x,0)\phi (x)=\mu _1\phi (x),\vspace{3pt%
} \\
\phi (x)\geq 0,\;\;\phi(x+L_i{\bf e}_i)=\phi(x), \, \forall 1\leq
i\leq n.
\end{array}
\right.  \label{2.3}
\end{equation}
If $\mu _1>0,$ then we say the zero solution is linearly unstable. In this
case, there exists a unique positive and $L$-periodic steady state solution $%
p(x)$ such that every positive solution $u(t,x)$ of (\ref{1}) with $L$%
-periodic initial function converges to $p(x)$ uniformly in $x\in
\mathbb{R}^n$ as $t\rightarrow \infty $, see \cite{B2005a}.

Now we can see that when $\mu _1>0,$ there are two non-negative
steady-states $p_1(x)=0$ and $p_2(x)=p(x),$ which has the same
mono-stability as in the Fisher-KPP equation. The heteroclinic connection
between $p_1(x)$ and $p_2(x)$ gives rise to a traveling wavefronts we are
seeking. In other words, by a traveling wave solution, we mean a particular
solution $u(t,x)=w(x,\xi ),\xi =x_1+ct$, which satisfies equation (\ref{1})
subject to the boundary condition
\begin{equation}
\lim_{s\rightarrow -\infty }w(x,s)=0,\;\lim_{s\rightarrow \infty
}w(x,\infty )=p(x),\text{ uniformly in }x\in \mathbb{R}^n.
\label{2.4}
\end{equation}
Substituting the relation $u(t,x)=w(x,\xi )$ with $\xi=x_1+ct$ into (\ref{1}), we obtain the
wave profile equation
\begin{equation}
\left\{
\begin{array}{l}
c\frac{\partial w}{\partial \xi }=\ \nabla \cdot (A(x)\nabla
w(x,x_1+ct))+f(x,w),\vspace{3pt} \\
w(x,-\infty )=0,\;w(x,\infty )=p(x),\;\;\;\text{uniformly in }x\in
\mathbb{R} ^n.
\end{array}
\right.  \label{2.5}
\end{equation}

Now we consider the case where $\mu _1>0,$ that is, the steady state
$0$ is unstable. The following result on the existence of spatially
periodic traveling wavefronts is from \cite {B2005,H2008}.

\noindent \textbf{Theorem A} (\cite{B2005,H2008}). Assume that $\mu
_1>0$ and $f$ satisfies (\ref{2.1}) and (\ref{2.2}). Then for any
$c\geq c^{*}$, there
exists a positive traveling wave $u(t,x)=w(x,x_1+ct)=w(x,\xi )$ satisfying (%
\ref{1}) and (\ref{2.4}), namely \eqref{2.5}, where the minimal speed $c^{*}$
is determined by
\begin{equation}
c^{*}=\min \{c:\;\exists \text{ }\lambda >0,\text{ such that }\mu _c(\lambda
)=0\}.  \label{2.5a}
\end{equation}
Here $\mu _c(\lambda )$ is the principal eigenvalue of the operator
\begin{eqnarray}
-L_{c,\lambda }\psi &=&-\nabla \cdot (A(x)\nabla \psi )-2\lambda \mathbf{%
e\cdot }\left( A(x)\nabla \psi \right)  \label{2.6} \\
&&-\left[ \lambda \nabla \cdot (A(x)\mathbf{e}^T)+\lambda ^2\mathbf{e\cdot }%
\left( A(x)\mathbf{e}^T\right) -\lambda c+f_u(x,0)\right] \psi ,  \nonumber
\end{eqnarray}
that is,
\[
-L_{c,\lambda }\psi =\mu _c(\lambda )\psi \ \mbox{ with } \,
\psi(x+L_i{\bf e}_i)=\psi(x), \, \forall  x\in \mathbb{R}^n,\, 1\leq
i\leq n.
\]

The derivation of $c^{*}$ in \cite{B2005,H2008} is based on the
technique developed by Aronson and Weinberger \cite{AW1,AW2} in
terms of a monotone-semiflow defined from equation (\ref{1}). Here
we can give an alternative, also easy to follow, explanation
directly from the wave profile equation (\ref{2.5}). This idea
originates from study of the KPP-Fisher equation and now is called
the technique of stability analysis of positive wavefronts near the
unstable steady state in their phase plane.

Due to biological or physical reason, we always expect a non-negative
wavefront, and as such, near the far-field $\xi =-\infty $ we assume $w\sim
v(x)e^{\lambda \xi }$ for some positive constant $\lambda $ and positive
periodic function $v(x),$ both of which are dependent on the wave speed $c.$
Substituting this relation into (\ref{2.5}) and ignoring the higher terms,
we then get
\begin{equation}
-\nabla \cdot \left( A(x)\nabla v\right) -2\lambda \mathbf{e\cdot }\left(
A(x)\nabla v\right) -\left[ \lambda \nabla \cdot (A(x)\mathbf{e}^T)+\lambda
^2\mathbf{e\cdot }\left( A(x)\mathbf{e}^T\right) -\lambda c+f_u(x,0)\right]
v=0,  \label{v}
\end{equation}
which implies that the operator $-L_{c,\lambda }$ defined in
(\ref{2.6}) should have zero as its principal eigenvalue. The
minimal speed $c$ is the smallest of such $c$ that (\ref{v}) holds
for $\lambda >0$. This agrees with the definition in (\ref{2.5a}).
Furthermore, in terms of $\lambda ,$ the operator $-L_{c,\lambda }$
is a polynomial of degree 2 and the result in \cite{H2008} shows
that the minimal speed $c^{*}$ is well-defined and has the following
properties.

\noindent \textbf{Theorem B} (\cite{H2008}). \textit{When $c>c^{*},$
there exists two constants $\lambda _1(c),\lambda
_2(c),\;\;0<\lambda _1(c)<\lambda _2(c)$ such that $\mu _c(\lambda
_1(c))=\mu _c\left( \lambda _2(c)\right) =0,\;\mu _c(\lambda )>0$
for any $\lambda \in (\lambda _1(c),\lambda _2(c))$; while in the
critical case $c=c^{*}$, we have $\mu _{c^*}(\lambda )\leq 0$, and
$\mu _{c^*}(\lambda ^{*})=0$ when $\lambda ^{*}=\lambda
_1(c^{*})=\lambda _2(c^{*})$.}

To present the exponential stability of wavefronts, we also impose a natural
condition on the positive steady-state $p_2=p(x).$ We assume that it is
linearly stable in the sense that the principal eigenvalue $\bar{\mu}_1$ of
the operator $L_1,$ defined by
\[
L_1\phi =\nabla \cdot (A(x)\nabla \phi )+f_u(x,p(x))\phi ,
\]
is negative, that is,
\begin{equation}
\left\{
\begin{array}{l}
\nabla \cdot (A(x)\nabla \phi )+f_u(x,p(x))\phi (x)=\bar{\mu}_1\phi (x),%
\vspace{3pt} \\
\phi (x)\geq 0,\;\;\phi(x+L_i{\bf e}_i)=\phi(x), \, \forall 1\leq
i\leq n
\end{array}
\right.  \label{2.7}
\end{equation}
hold for some negative constant $\bar{\mu}_1$ and non-zero function
$\phi$.

The subsequent result shows that the condition $\bar{\mu}_1<0$ can
be easily realized.

\begin{lem}
Assume that
\begin{enumerate}
\item[(C)] $f_u(x,u)\leq \frac{f(x,u)}{u},\, \forall u>0, \, x\in \mathbb{R}^n$, and
there exists $x_0\in \mathbb{R}^n$ such that $f_u(x_0,u)<
\frac{f(x_0,u)}{u}, \, \forall u>0$.
\end{enumerate}
Then any possible positive $L$-periodic steady state $p(x)$ of
(\ref{1}) is linearly stable.
\end{lem}

\begin{proof}
Let $g(x,u):=f(x,u)/u,\, \forall u>0$, and let $\mu (m(x))$ be the
principal eigenvalue of the operator $L_1$ with $f_u(x,p(x))$
replaced by a continuous and $L$-periodic weight function $m(x)$.
Since $p(x)=p(x+L)$ satisfies the elliptic equation and only the
principal eigenvalue has positive eigenfunction, it is easy to see
that $\mu (g(x,p(x))=0$. Since $f_u(x,p(x))\leq g(x,p(x))$  and
$f_u(x_0,p(x_0))<g(x_0,p(x_0))$, the monotonicity of the principal
eigenvalue with respect to the weight function implies that $\mu
(f_u(x,p(x)))<\mu (g(x,p(x)))=0$. Thus, $p(x)$ is linearly stable.
\end{proof}

For any given wavefront $w(x,\xi)$,  $\xi=x_1+ct$ with speed $c\geq
c^{*}$, we define a weight function
\begin{equation}
W(\xi )=\begin{cases} e^{-\lambda (\xi -\xi _0)}, & \xi \leq \xi _0, \ \ \
\lambda \in \left[ \lambda _1(c),\lambda _2(c)\right] ,\vspace{3pt}\; \\
1,&\xi > \xi _0, \end{cases}  \label{W}
\end{equation}
where $\xi _0$ is chosen in such a way that the wavefront $w(x,\xi ),$ $\xi
\geq \xi _0,$ is very close to the positive steady-state $p_2=p(x)$, that
is, $\xi _0$ is chosen so that
\begin{equation}
|w(x,\xi )-p(x)|<\bar{\varepsilon} ,\;\;\forall \xi \geq \xi _0,
\label{2.11a}
\end{equation}
uniformly in $x\in \mathbb{R}^n$ for any given sufficient small
constant $\bar{\varepsilon}
>0.$

Now we are in a position to prove the following stability result.

\begin{thm}
\label{thm2.1} Assume that \eqref{2.1} and \eqref{2.2} hold,  $\mu
_1>0$, and $\bar{\mu}_1<0$.
Let $u(t,x)$ be the  solution (\ref{1}) with the initial condition $%
u_0(x)$ satisfying $0\leq u_0(x)\leq p(x)$, and $w(x,x_1+ct)$ be a
spatially periodic traveling wavefront with $c\geq c^{*}.$ If the
initial data $u_0(x),$ having at most finite points of
discontinuity, satisfies
\begin{equation}
W(x_1)(u_0(x)-w(x,x_1))\in L^1(\mathbb{R}^n\mathbf{),}  \label{2.9}
\end{equation}
then the following statements are valid:

(1) For $\xi =x_1+ct\leq \xi _0,$ it holds
\begin{equation}
|u(t,x)-w(x,x_1+ct)|<C(1+t)^{-\frac n2}e^{-\mu _c(\lambda )t};  \label{2.10}
\end{equation}

(2) For $\xi =x_1+ct\geq \xi _0,$ it holds
\begin{equation}
|u(t,x)-w(x,x_1+ct)|<C(1+t)^{-\frac n2}e^{-\min \{\mu _c(\lambda ),\frac{-%
\bar{\mu}_1}2\}t}.  \label{2.11}
\end{equation}
Here $\lambda $, as mentioned in the weight function, is any number in the
interval $\left[ \lambda _1(c),\lambda _2(c)\right] .$
\end{thm}

\begin{proof}
For the given initial function $u_0(x)$, we define
\begin{equation}
\left\{
\begin{array}{l}
U_0^{+}(x)=\max \{u_0(x),\;w(x,x_1)\},\;x\in \mathbb{R}^n\mathbf{,\vspace{3pt%
}} \\
U_0^{-}(x)=\min \{u_0(x),\;w(x,x_1)\},\;x\in \mathbb{R}^n,
\end{array}
\right.  \label{2.14}
\end{equation}
where $w(x,x_1)$ is the initial function of the traveling wave
solution $w(x,x_1+ct)$.

Let $U^{\pm }(t,x)$ be the solutions of (\ref{1}) with the initial data $%
U_0^{\pm }(x),$ respectively, that is,
\begin{equation}
\left\{
\begin{array}{l}
\dfrac{\partial U^{\pm }}{\partial t}=\nabla \cdot \left( A(x)\nabla U^{\pm
}\right) +f(x,U^{\pm }),\vspace{4pt} \\
U^{\pm }(0,x)=U_0^{\pm }(x),\;x\in \mathbb{R}^n\mathbf{.}
\end{array}
\right.  \label{2.14a}
\end{equation}
Since $0\leq U_0^{-}\leq w(x,x_1)\leq U_0^{+}\leq p(x),$ the comparison
principle implies that
\[
0\leq U^{-}(t,x)\leq w(x,x_1+ct)\leq U^{+}(t,x)\leq p(x),
\]
and
\[
0\leq U^{-}(t,x)\leq u(t,x)\leq U^{+}(t,x)\leq p(x).
\]
Obviously, we have
\[
U^{-}(t,x)-w(x,x_1+ct)\leq u(t,x)-w(x,x_1+ct)\leq U^{+}(t,x)-w(x,x_1+ct),
\]
or
\[
|u(t,x)-w(x,x_1+ct)|\leq \max \{|U^{+}(t,x)-w(x,x_1+ct)|,\
|U^{-}(t,x)-w(x,x_1+ct)|\}.
\]
As such, to prove (\ref{2.10}) and (\ref{2.11}), we proceed to show
\begin{equation}
|U^{\pm }(t,x)-w(x,x_1+ct)|\leq C(1+t)^{-\frac n2}e^{-\mu _c(\lambda )t},\
\text{ for }\ \xi \leq \xi _0  \label{2.15}
\end{equation}
and
\begin{equation}
|U^{\pm }(t,x)-w(x,x_1+ct)|\leq C(1+t)^{-\frac n2}e^{-\min \{\mu _c(\lambda
),\frac{-\bar{\mu}_1}2\}t},\;\;\text{for }\;\xi >\xi _0.  \label{2.16}
\end{equation}
In what follows, we give a detailed proof of (\ref{2.15}) and (\ref{2.16})
for $U^{+}(t,x)$. The similar arguments work for $U^{-}(t,x)$, too, and hence, we
omit them.

Recall that both $U^{+}(t,x)$ and $w(x,x_1+ct)$ satisfy the first equation
of (\ref{2.14a}). Note that $\xi =x_1+ct.$ Thus we define
\[
V(t,x)=U^{+}(t,x)-w(x,x_1+ct).
\]
It is easy to see that $V(t,x)$ is a nonnegative solution of the
following equation
\[
\begin{cases}
V_t=\nabla \cdot \Big( A \nabla  V \Big)  +f(x,w+V)-f(x,w)
\vspace{4pt} \\
V(0,x )=U_0^{+}(x )-w(x ,x_1 )\geq 0,
\end{cases}
\]
with the initial data $V(0,x)$ satisfying
\begin{equation}
W(x_1)V(0,x)=W(x_1)\left( U_0^{+}(x)-w(x,x_1)\right) \in
L^1(\mathbb{R}^n), \label{2.17}
\end{equation}
where $W(x_1)$ is the weight function defined in (\ref{W}).

Since the non-linear function $f$ is sub-linear and satisfies (\ref{2.1a}),
it follows that
\begin{equation}
f(x,w+V)-f(x,w)\leq f_u (x,0)V.  \label{Q?}
\end{equation}
We then have an inequality
\begin{equation}
\begin{cases} V_t\leq \nabla \cdot \Big( A \nabla V \Big) +f_u (x,0)V,\vspace{4pt} \\ V(0,x )=U_0^{+}(x )-w(x ,x_1 )\geq 0. \end{cases}
\label{2.18}
\end{equation}

Let $v(x),$ dependent on $c$ and $\lambda ,$ be the eigenfunction of
operator $-L_{c,\lambda }$ with the principal eigenvalue $\mu _c(\lambda ).$
From Theorem B, when $c\geq c^{*},$ $\mu _c(\lambda )$ is real and $v(x)>0$
is well-defined. By a crucial wave transformation
\begin{equation}
V(t,x)=v(x)e^{\lambda (x_1+ct -\xi _0)}\bar{V}(t,x)e^{-\mu
_c(\lambda )t},\;\;\lambda \in [\lambda _1(c),\lambda _2(c)],
\label{2.19}
\end{equation}
it then follows from (\ref{2.18}) that $\bar{V}(t,x)$ satisfies
\begin{equation}
\begin{cases} \bar{V}_t\leq \nabla \cdot \Big( A \nabla \bar{V} \Big)
+g(t,x) \cdot \nabla \bar{V},\vspace{4pt} \\
\bar{V}(0,x)=v^{-1}(x)e^{-\lambda (x_1 -\xi _0)}V(0,x)\geq 0,
\end{cases}  \label{2.20}
\end{equation}
where
\begin{equation}
g=2A\left( \frac{\nabla v}v+\lambda \mathbf{e}\right) .  \label{2.21}
\end{equation}
To get (\ref{2.20}), we have made use of the following formulas
\[
\nabla (uv)=u\nabla v+v\nabla u,\;\;\;\;\nabla \cdot (vA\nabla u)=v\nabla
\cdot (A\nabla u)+\nabla v\cdot (A\nabla u).
\]
It is easy to see that $\bar{V}$ in (\ref{2.20}) can be estimated by
\begin{equation}
\bar{V}(t,x)\leq \int_{\mathbb{R}^n}G(t,x-y)\bar{V}(0,y)dy,\quad
\forall x\in \mathbb{R}^n,\, t\geq 0, \label{2.22}
\end{equation}
where $G(t,x-y)$ is the fundamental solution (i.e., the Green function) so
that $\int G(t,x-y)u_0(y)dy$ satisfies the following partial differential
equation
\[
u_t=\nabla \cdot (A\nabla u)-g\cdot \nabla u,\;\;\;\;u(0,x)=u_0(x),
\]
with
\[
\int_{\mathbf{R}}G(t,y)dy=1.
\]
Furthermore, since $A$ is spatially periodic and uniformly elliptic,
it follows from \cite[(6.12) in page 24]{FR} that
\begin{equation}
|G(t,x-y)|\leq \frac C{t^{\frac n2}}e^{-\frac{|x-y|^2}{4\alpha _1t}}
\label{2.23}
\end{equation}
for some $\alpha _1>0$, where $n$ is the spatial dimension. By
(\ref{2.17}), we understand that $\bar{V}(0,x)\in
L^1(\mathbb{R}^n).$ From (\ref{2.22}), together with (\ref{2.23}),
we can obtain
\[
\bar{V}(t,x)\leq Ct^{-\frac
n2}\int_{\mathbb{R}^n}|\bar{V}(0,y)|dy,\quad \forall x\in
\mathbb{R}^n, \, t\geq 0.
\]
Since the solution $\bar{V}$ for $t$ near zero has no singularity, the term $%
t^{-\frac n2\text{ }}$ can be replaced by $(1+t)^{-n/2}$. This technique is
frequently used in the heat equation. Returning to (\ref{2.19}), we then
have
\begin{eqnarray}
V(t,x) &=&|U^{+}(t,x)-w(x,x_1+ct)|  \nonumber \\
&=&v(x)e^{\lambda (x_1+ct -\xi _0)}\bar{V}(t,x)e^{-\mu _c(\lambda
)t} \nonumber
\\
&\leq &Ce^{\lambda (x_1+ct -\xi _0)}(1+t)^{-\frac{n}{2}}e^{-\mu
_c(\lambda )t} \label{vv}
\end{eqnarray}
for all $x\in \mathbb{R}^n, \, t\geq 0$. When $\xi=x_1+ct \leq \xi _0$, we have $%
e^{\lambda (\xi -\xi _0)}\leq 1$. This proves (\ref{2.10}).

Next for $\xi=x_1+ct >\xi _0,$ we know that
\[
w(x,\xi )\leq w(x,\xi)+V\leq p(x).
\]
Since we have chosen $\xi _0$ so that $w(x,\xi )$ is very close to
$p(x)$ uniformly in $x\in \mathbb{R}^n$ when $\xi >\xi _0$ (for
precision, see (\ref{2.11a})), as such for $\xi =x_1+ct >\xi _0 $,
by the regularity of $f$, we have
\[
f(x,w+V)-f(x,w)\leq (f_u(x,p(x))+\varepsilon )V
\]
and
\[
\begin{cases}
V_t\leq \nabla \cdot (A\nabla V)
  +(f_u (x,p(x))+\varepsilon
 )V,
\vspace{4pt} \\
V(0,x)=U_0^{+}(x )-w(x,x_1)\geq 0
\end{cases}
\]
for some small constant $\varepsilon $ satisfying $0<\varepsilon <-\bar{\mu}%
_1$, where $\bar{\mu}_1$ is defined in (\ref{2.7}) with
eigenfunction $\phi (x)>0$. Under the transformation
\begin{equation} \label{moz1}
V(t,x)=\bar{V}(t,x)\phi (x)e^{(\bar{\mu}_1+\varepsilon )t},
\end{equation}
it then follows that for all $\xi=x_1+ct >\xi _0$, $\bar{V}(t,x)$
satisfies
\begin{equation}
\bar{V}_t\leq \nabla \cdot (A\nabla \bar{V})+2\left( A\frac{\nabla \phi }%
\phi \right) \cdot \nabla \bar{V},\ \ \ \ \bar{V}(0,x)=\phi ^{-1}(x)V(0,x),
\label{2.24}
\end{equation}
with
\[
\bar{V}(t,x)\leq V(t,x)\phi ^{-1}(x)e^{-(\bar{\mu}_1+\varepsilon )t}
\]
on the line $x_1+ct=\xi_0$.  Since on the boundary line $\xi
=x_1+ct=\xi _0,$ (\ref{vv}) yields
\begin{equation}\label{moz2}
V(t,x)\leq C(1+t)^{-\frac n2}e^{-\mu _c(\lambda )t},
\end{equation}
which implies that on the line $x_1+ct=\xi _0,$
\[
\bar{V}(t,x)\leq C(1+t)^{-\frac n2}e^{-\mu _c(\lambda )t}e^{\mu
_2t},
\]
where $\mu _2=-\bar{\mu}_1-\varepsilon >0$.  Let
\[
\mu _{\min }=\min \{\mu _c,\frac{-\bar{\mu}_1}2\}.
\]
Thus, applying \eqref{moz2} to \eqref{moz1}, we obtain that
 on the line $x_1+ct=\xi _0,$
\[
\bar{V}(t,x)\leq C(1+t)^{-\frac n2}e^{-\mu _{\min }t}e^{\mu _2t}.
\]
Obviously, $C(1+t)^{-\frac n2}e^{-\mu _{\min }t}e^{\mu _2t}$ is an upper
solution of the equation
\[
u_t=\nabla \cdot (A\nabla u)+(2A\frac{\nabla \phi }\phi )\cdot
\nabla u,\quad t\geq 0
\]
for some constant $C$. By the parabolic comparison principle, it
then follows that
\[
V(t,x) =\bar{V}(t,x)\phi (x)e^{(\bar{\mu}_1+\varepsilon )t} \leq
C(1+t)^{-\frac n2}e^{-\mu _{\min }t},\, \,  \forall x\in
\mathbb{R}^n,\, t\geq 0.
\]
This completes the proof.
\end{proof}

The following result is a straightforward consequence of Theorem
\ref{thm2.1}.

\begin{thm}
\label{thm2.2} Assume that \eqref{2.1} and \eqref{2.2} hold,
$\mu_1>0$, and $\bar{\mu}_1<0$.
Let $u(t,x)$ be the solution (\ref{1}) with the initial condition $%
u_0(x),\;0\leq u_0(x)\leq p(x),$ and $w(x,x_1+ct)$ be a traveling
wavefront with $c\geq c^{*}.$ If the initial data $u_0(x),$ having
at most finite points of discontinuity, satisfies
\[
W(x_1)(u_0(x)-w(x,x_1))\in L^1(\mathbb{R}^n\mathbf{)},
\]
then the following statements are valid:

(1) For any $c>c^{*}$ and $\lambda $ $\in (\lambda _1(c),\lambda _2(c)),$
the wavefront $w(x,x+ct)$ is exponentially stable in the sense that
\begin{equation}
\sup_{x\in \mathbb{R}^n}|u(t,x)-w(x,x_1+ct)|\leq Ce^{-\mu
_3t},\;\;\quad \text{for some }\,\,\mu _3>0;  \label{2.12}
\end{equation}

(2) For $c=c^{*}$, or $c>c^{*}$ with $\lambda =\lambda _1(c),$ the wavefront
is algebraically stable in the sense that
\begin{equation}
\sup_{x\in \mathbb{R}^n}|u(t,x)-w(x,x_1+ct)|\leq C(1+t)^{-\frac n2}.
\label{2.13}
\end{equation}
\end{thm}

\begin{proof}
For any given $c>c^{*}$ and $\lambda \in (\lambda _1(c),\lambda
_2(c)),$
it follows from Theorem B that $\mu _c(\lambda )>0.$ Let $\mu _3=\min \{$ $%
\mu _c(\lambda ),\mu _2\}.$ Then we have (\ref{2.12}) from (\ref{2.10}) and (%
\ref{2.11}). When $c=c^{*},$ or $c>c^{*}$ with $\lambda =\lambda
_1(c),$ we then get $\mu _c(\lambda )=0,$ and a combination of
(\ref{2.10}) and (\ref {2.11}) yields (\ref{2.13}).
\end{proof}

\section{An extension to the time-periodic case}

We have already studied the system which is periodic in space. In this
section, we extend our results to the case when the system is also periodic
in time.

We consider the following parabolic equation
\begin{equation}
u_t-\nabla \cdot (A(t,x)\nabla u)=f(t,x,u),\;\;x\in \mathbb{R}^n,
\label{5.1}
\end{equation}
where the matrix $A$ and the nonlinear term $f$ are periodic both in $t$ and
$x,$ that is,
\[
A(t+T,x)=A(t,x),\;\;f(t+T,x,u)=f(t,x,u),\;
\]
and
\[
A(t,x+L_i{\bf e}_i)=A(t,x),\;f(t,x+L_i{\bf e}_i,u)=f(t,x,u),\, \,
\forall 1\leq i\leq n
\]
for some positive real number $T>0$ and $n$-dimensional vector $%
L=(L_1,L_2,...,L_n)$.

As in \cite[page 364-366]{H2008}, we assume that  the uniformly
elliptic symmetric matrix $A(t,x)$ is of class $C_{t,x}^{1,\alpha
/2;1,\alpha }$ and the function $f(t,x,u)$ is continuous, of class
$C^{0,\alpha /2;0,\alpha }$ with respect to $(t,x)$ locally
uniformly in $u\in \mathbb{R}$. We also assume that $\frac{\partial
f}{\partial u}$ exists and is continuous in $u$ uniformly for
$x\in\mathbb{R}^n$ and $t \in \mathbb{R}$ . The existence of
traveling waves $w(t,x,x_1+ct)$ was studied in \cite[page
364-366]{H2008} and \cite{nadin} with the sub-linearity assumption
on the reaction
term $f.$ Assume that there are two periodic states $p_1(t,x)=0$ and $%
p_2(t,x)>0$ of (\ref{5.1}). The system is called mono-stable if the state $%
p_1(t,x)$ is unstable and the other state $p_2(t,x)$ is stable. In what
follows, we assume that $\mu _1>0$ and $\bar{\mu}_1<0,$ where $\mu _1$ is
the principal eigenvalue of the linearized operator around $0$:
\[
\varphi (t,x)\rightarrow \varphi _t-\nabla \cdot (A(t,x)\nabla \varphi
)-f_u(t,x,0)\varphi ,\;\;\varphi (t,x)\text{ is periodic in }t\text{ and }x,
\]
and $\bar{\mu}_1$ is the the principal eigenvalue of the linearized operator
around $p_2(t,x)$:
\[
\varphi (t,x)\rightarrow \varphi _t-\nabla \cdot (A(t,x)\nabla \varphi
)-f_u(t,x,p_2)\varphi ,\;\;\varphi (t,x)\text{ is periodic in }t\text{ and }%
x.
\]
Traveling wavefronts are solutions of (\ref{5.1}) in the form $u(t,x)=w(t,x,%
\mathbf{e}\cdot x+ct)=w(t,x,x_1+ct),$ where $\mathbf{e}=(1,0,...,0),$ which
is periodic in the first two variables $t$ and $x$ and satisfies
\[
\lim_{\xi \rightarrow -\infty }w(t,x,\xi )=p_1(t,x)=0,\;\;\lim_{\xi
\rightarrow \infty }w(t,x,\xi )=p_2(t,x)>0
\]
uniformly for all $t$ and $x.$ The minimal speed $c^{*}$ is defined by
\[
c^{*}=\min \{c:\;\exists \text{ }\lambda \text{ so that }\mu _c(\lambda
)=0\},
\]
where $\mu (\lambda )$ is the principal eigenvalue of the operator
\[
v(t,x)\rightarrow v_t+c\lambda v-\nabla \cdot (A(t,x)\nabla v)-2\lambda
\mathbf{e\cdot }\left( A\nabla v\right) -\lambda \nabla \cdot (A\mathbf{e}%
)-\lambda ^2\mathbf{e\cdot }\left( A\mathbf{e}\right) -f_u(t,x,0)v,\;
\]
that is, there exists a periodic function $v(t,x)>0$ such that
\[
v_t+c\lambda v-\nabla \cdot (A(t,x)\nabla v)-2\lambda \mathbf{e\cdot }\left(
A\nabla v\right) -\lambda \nabla \cdot (A\mathbf{e})-\lambda ^2\mathbf{%
e\cdot }\left( A\mathbf{e}\right) \ -f_u(t,x,0)v=\mu _c(\lambda )v.
\]
It then follows that when $c>c^{*},\;$there exist $0<\lambda _1(c)<\;\lambda
_2(c)$ such that $\mu _c(\lambda _1)=\mu _c(\lambda _2)=0$ and $\mu
_c(\lambda )>0$ for $\lambda \in (\lambda _1,\lambda _2).$ In particular,
when $c=c^{*},$ we have $\lambda _1(c^{*})=\;\lambda _2(c^{*})$ and $\mu
_c(\lambda )\leq 0$ for all $\lambda .$

The weight function is defined as in (\ref{W}). The inequality (\ref{2.18})
now becomes
\[
\left\{
\begin{array}{l}
V_t\leq \nabla \cdot \left( A\nabla V\right) +f_u^{\prime }(t,x,0)V,\vspace{%
4pt} \\
V(0,x)=U_0^{+}(x)-w(0,x,x_1)\geq 0.
\end{array}
\right.
\]
Correspondingly, the transformation (\ref{2.19}) is replaced by
\[
V(t,\xi )=v(t,x)e^{\lambda (x_1+ct -\xi _0)}\bar{V}(t,x)e^{-\mu
_c(\lambda )t},\;\;\lambda \in [\lambda _1(c),\lambda _2(c)].
\]
All the remaining arguments in section 2 still work. Therefore, we have the
following results.

\begin{thm}
Let $u(t,x)$ be a solution (\ref{5.1}) with initial condition $%
u_0(x),\;0\leq u_0(x)\leq p_2(t,x),$ and $w(t,x,x_1+ct)$ be a traveling
wavefront with wave speed $c\geq c^{*}.$ If the initial data $u_0(x),$
having at most finite points of discontinuity, satisfies
\begin{equation*}
W(x_1)(u_0(x)-w(0,x,x_1))\in L^1(\mathbb{R}^n\mathbf{),}
\label{5.3}
\end{equation*}
then the following statements are valid:

(1) For $\xi =x_1+ct\leq \xi _0$, it holds
\begin{equation*}
|u(t,x)-w(t,x,\xi )|<C(1+t)^{-\frac n2}e^{-\mu _c(\lambda )t};  \label{5.4}
\end{equation*}

(2) For $\xi =x_1+ct\geq \xi _0$, it holds
\begin{equation*}
|u(t,x)-w(t,x,\xi )|<C(1+t)^{-\frac n2}e^{-\min \{\mu _c(\lambda ),\frac{-%
\bar{\mu}_1}2\}t}.  \label{5.5}
\end{equation*}
Here $\lambda $ is any number in the interval $\left[ \lambda _1(c),\lambda
_2(c)\right] .$
\end{thm}

\begin{thm}
Let $u(t,x)$ be a solution (\ref{5.1}) with initial condition $%
u_0(x),\;0\leq u_0(x)\leq p_2(t,x),$ and $w(t,x,x_1+ct)$ be a
traveling wavefront with wave speed $c\geq c^{*}.$ If the initial
data $u_0(x),$ having at most finite points of discontinuity,
satisfies
\[
W(x_1)(u_0(x)-w(0,x,x_1))\in L^1(\mathbb{R}^n\mathbf{),}
\]
then the following statements are valid:

(1) For any $c>c^{*}$ and $\lambda \in (\lambda _1(c),\lambda _2(c)),$ the
wavefront $w(t,x,x_1+ct)$ is exponentially stable in the sense that
\begin{equation*}
\sup_{x\in \mathbb{R}^n}|u(t,x)-w(t,x,x_1+ct)|\leq Ce^{-\mu
_3t},\;\;\quad \text{for some}\,\,\,\mu _3>0;  \label{5.6}
\end{equation*}

(2) For $c=c^{*}$, or $c>c^{*}$ with $\lambda =\lambda _1(c),$ the wavefront
is algebraically stable in the sense that
\begin{equation*}
\sup_{x\in \mathbb{R}^n}|u(t,x)-w(t,x,x_1+ct)|\leq C(1+t)^{-\frac
n2}. \label{5.7}
\end{equation*}
\end{thm}

\section{A concluding remark}

The stability analysis above, together with asymptotic behavior of
spatially periodic wavefronts $w(x,x_1+ct)$ near the far field, can
be directly used to obtain the uniqueness of the wave with a given
speed. This then may provide a new and insightful approach that is
different from the technique in \cite{H2}.

To illustrate this idea, we first recall a result from \cite[Theorem 1.3]
{H2008}.

\noindent \textbf{Theorem C} (\cite{H2008}). \textit{If }$c>c^{*},$\textit{\
then there exists a constant }$B>0$\textit{\ such that }
\begin{equation}
w(x,\xi )\sim Be^{\lambda _1(c)\xi }v(x),\;\;\xi \rightarrow -\infty
\label{3.1}
\end{equation}
\textit{where }$v(x)$\textit{\ is the eigenfunction defined in (\ref{v}). If
}$c=c^{*},$\textit{\ then there exists a constant }$B>0$\textit{\ such that }
\begin{equation}
w(x,\xi )\sim B|\xi |e^{\lambda _1(c^{*})\xi }v(x),\;\;\xi \rightarrow
-\infty  \label{3.2}
\end{equation}
\textit{where }$\lambda _1(c^{*})=\lambda _2(c^{*})=\lambda ^{*}$\textit{\
and }$v(x)$\textit{\ is the eigenfunction in (\ref{v}) with }$c=c^{*}.$

Indeed, we can further claim that if $c>c^{*},$ then the wave profile $%
w(x,\xi ) $ has the following asymptotical behavior:

\begin{equation}
w(x,\xi )=Be^{\lambda _1(c)\xi }v(x)+O(e^{\alpha \xi }),\;\;\;B>0,\;\alpha
>\lambda _1(c),\;\;\;\xi \rightarrow -\infty .  \label{3.3}
\end{equation}

To see this, we may assume from (\ref{3.2}) that $w(x,\xi )=Be^{\lambda
_1(c)\xi }v(x)(1+$ $\bar{w}(x,\xi )).$ Substituting this into (\ref{2.5}),
we can get an equation for $\bar{w}(x,\xi ).$ For this new equation, after a
simple asymptotic analysis near $\xi =-\infty ,$ we can find that $\bar{w}%
(x,\xi )$ decays to zero exponentially as $\xi \rightarrow -\infty
.$ This confirms our claim. Now for any given $c>c^{*},$ suppose
that there are two spatially periodic wave profiles $w_1(x,\xi )$
and $w_2(x,\xi )$ of (\ref{1}). A shift of distance $s$ in $w_2$
then enables $w_2(x,\xi +s)$ to share the same leading term in
(\ref{3.3}) as that of $w_1.$ Let $\lambda $ in the weight function
be $\lambda _1(c).$ Then when $t=0,$ we have $e^{-\lambda
_1(c)(x_1-\xi _0)}(w_1(x,x_1)-w_2(x,x_1+s))\in L^1(\mathbb{R}^n)$.
By Theorem \ref{thm2.2}, it follows that
\[
\lim_{t\to \infty}\sup_{x\in
\mathbb{R}^n}\left|w_1(x,x_1+ct)-w_2(x,x_1+ct +s)\right|=0.
\]
Using a change of variable $\xi=x_1+ct$, we then see that
\begin{equation}\label{unique}
\lim_{t\to \infty}\sup_{(\xi,x_2,\cdots,x_n)\in
\mathbb{R}^n}\left|w_1((\xi-ct,x_2,\cdots,
x_n),\xi)-w_2((\xi-ct,x_2,\cdots, x_n),\xi+s)\right|=0.
\end{equation}
By the $L_1$-periodicity of $w_1((x_1,x_2,\cdots,x_n),\xi)$ and
$w_2((x_1,x_2,\cdots, x_n),\xi+s)$ with respect to $x_1\in
\mathbb{R}$, it easily follows from (\ref{unique}) that
\[
w_1(x,\xi)=w_2(x,\xi +s),\,
\, \forall x=(x_1,x_2,\cdots, x_n)\in \mathbb{R}^n,\, \xi\in
\mathbb{R}.
\]
This shows that for any given $c>c^{*}$, the wave profile $w(x,\xi
)$ is unique up to translation in $\xi $. In the case where
$c=c^{*},$ the same idea applies, but more delicate details are
needed.

\

\noindent {\bf Acknowledgements.}  We are very grateful to two
anonymous referees for careful reading and helpful suggestions which
led to an improvement of our original manuscript.

\end{document}